\documentclass{amsart}
\usepackage{amscd,amssymb,url}
\usepackage[matrix,arrow]{xy}
\theoremstyle{plain}
\newtheorem{Thm}{Theorem}
\newtheorem*{MainThm}{Main Theorem}
\newtheorem{Lem}[Thm]{Lemma}
\newtheorem{Prop}[Thm]{Proposition}

\DeclareMathOperator{\HA}{\mathcal H}
\DeclareMathOperator{\modcat}{mod}
\DeclareMathOperator{\indcat}{ind}
\DeclareMathOperator{\dimv}{\underline\dim}
\DeclareMathOperator{\Hom}{Hom}
\DeclareMathOperator{\End}{End}
\DeclareMathOperator{\Ext}{Ext}
\DeclareMathOperator{\Aut}{Aut}
\begin{document}

\title{Hall Polynomials for Affine Quivers}
\author{Andrew Hubery}
\address{Universit\"at Paderborn\\Germany}
\email{hubery@math.upb.de}

\begin{abstract}
We use Green's comultiplication formula to prove that Hall polynomials exist for all Dynkin and affine quivers. For Dynkin and cyclic quivers this approach provides a new and simple proof of the existence of Hall polynomials. For non-cyclic affine quivers these polynomials are defined with respect to the decomposition classes of Bongartz and Dudek, a generalisation of the Segre classes for square matrices. 
\end{abstract}

\maketitle

In \cite{Ringel1}, Ringel showed how to construct an associative algebra from the category of finite modules over a finitary ring, and whose multiplication encodes the possible extensions of modules. More precisely, one forms the free abelian group with basis the isomorphism classes of modules and defines a multiplication by taking as structure constants the Hall numbers
$$F_{MN}^X:=|\{U\leq X:U\cong N, X/U\cong M\}|.$$
Alternatively, we can write $F_{MN}^X=P_{MN}^X/a_Ma_N$, where
$$P_{MN}^X:=|\{(f,g):0\to N\xrightarrow{f} X\xrightarrow{g} M\to 0 \textrm{ exact}\}| \quad\textrm{and}\quad a_M:=|\Aut(M)|.$$
In the special case of finite length modules over a discrete valuation ring with finite residue field, one recovers the classical Hall algebra \cite{Macdonald}.

Interesting examples of such Hall algebras arise when one considers representations of a quiver over a finite field. Green showed in \cite{Green} that the subalgebra generated by the simple modules, the so-called composition algebra, is a specialisation of Lusztig's form for the associated quantum group. (More precisely, one must first twist the multiplication in the Hall algebra using the Euler characteristic of the module category.) In proving this result, Green first showed that the Hall algebra is naturally a self-dual Hopf algebra.

If one starts with a Dynkin quiver, then the isomorphism classes of indecomposable modules are in bijection with the set of positive roots $\Phi_+$ of the corresponding semisimple complex Lie algebra, the bijection being given by the dimension vector \cite{Gabriel}. Thus the set of isomorphism classes of modules is in bijection with the set of functions $\xi\colon\Phi_+\to\mathbb N_0$, and as such is independent of the field.

In this setting, Ringel proved in \cite{Ringel2} that the Hall numbers are given by universal polynomials; that is, given three functions $\mu,\nu,\xi\colon\Phi_+\to\mathbb N_0$, there exists a polynomial $F_{\mu\nu}^\xi\in\mathbb Z[T]$ (depending only on $Q$) such that, for any finite field $k$ with $|k|=q$ and any representations $M$, $N$ and $X$ belonging to the isomorphism classes $\mu$, $\nu$ and $\xi$ respectively, we have
$$F_{MN}^X = F_{\mu\nu}^\xi(q).$$
The proof uses associativity and induction on the Auslander-Reiten quiver, reducing to the case when $M$ is isotypic.
A similar result holds when $Q$ is an oriented cycle and we consider only nilpotent modules. The indecomposable nilpotent representations are all uniserial, so determined by their simple top and Leowy length. In fact, if $Q$ is the Jordan quiver, consisting of a single vertex and a single loop, then the Ringel-Hall algebra is precisely the classical Hall algebra, so isomorphic to Macdonald's ring of symmetric functions. It is well-known that Hall polynomials exist in this context \cite{Macdonald}.
A natural question to ask, therefore, is whether Hall polynomials exist for more general quivers. Certain results along this line are clear, for example if all three modules are preprojective or preinjective. Also, some Hall numbers have been calculated for modules over the Kronecker quiver \cite{Szanto} (or equivalently for coherent sheaves over the projective line \cite{BK}) and we see that there is again `polynomial behaviour'. For example, let us denote the indecomposable preprojectives by $P_r$ and the indecomposable preinjectives by $I_r$, for $r\geq0$. If $R$ is a regular module of dimension vector $(n+1)\delta$ and containing at most one indecomposable summand from each tube, then
$$P_{P_mI_{n-m}}^R = (q-1)a_R = P_{RP_m}^{P_{m+n+1}},$$
so that
$$F_{P_mI_{n-m}}^R = a_R/(q-1) \quad\textrm{and}\quad F_{RP_m}^{P_{m+n+1}}=1.$$
Thus, although the regular modules $R$ depend on the field, the Hall numbers depend only on the number of automorphisms of $R$.

In general, the isomorphism classes of indecomposable $kQ$-modules are no longer combinatorially defined --- they depend intrinsically on the base field $k$. Some care has therefore to be taken over the definition (and meaning) of Hall polynomials.
The existence of Hall polynomials has gained importance recently by the relevance of quiver Grassmannians and the numbers
$$\big|\mathrm{Gr}\textstyle{\binom{X}{\underline e}}\big|:=\displaystyle{\sum_{\substack{[M],[N]\\\dimv N=\underline e}}F_{MN}^X}$$
to cluster algebras. In \cite{CR} Caldero and Reineke show for affine quivers that these numbers are given by universal polynomials. Also, as shown in \cite{Reineke}, the existence of universal polynomials implies certain conditions on the eigenvalues of the Frobenius morphism on $l$-adic cohomology for the corresponding varieties.
The aim of this article is to show how the comultiplication, or Green's Formula, can be used to prove the existence of Hall polynomials. We remark that Ringel's proof cannot be extended, since the associativity formula alone doesn't reduce the difficulty of the problem --- the middle term remains unchanged. The advantage of Green's Formula is that it reduces the dimension vector of the middle term and hence allows one to apply induction. In this way we can reduce to a situation where the result is clear: for Dynkin quivers we reduce to the case when the middle term is simple; for nilpotent representations of a cyclic quiver we reduce to the case when the middle term is indecomposable; for general affine quivers, we reduce to the case when either all three representations are regular, or else the middle term is regular and the end terms are indecomposable. In this latter case, we can use associativity and a result of Schofield on exceptional modules \cite{Ringel5} to simplify further to the case when $M$ is simple preinjective.

After recalling the necessary theory, we apply our reductions in the special cases of Dynkin quivers and nilpotent representations of a cyclic quiver, cases which are of course of particular interest. This gives a new proof of the existence of Hall polynomials in these cases, and the proofs offered here are short and elementary.

We then extend this result to all finite dimensional representations of a cyclic quiver, where we see that Hall polynomials exist with respect to the Segre classes. More precisely, given three Segre symbols $\alpha$, $\beta$ and $\gamma$, there is a universal polynomial $F_{\alpha\beta}^\gamma$ such that, for any finite field $k$ with $|k|=q$,
$$\sum_{\substack{A\in\mathcal S(\alpha,k)\\B\in\mathcal S(\beta,k)}}F_{AB}^C = F_{\alpha\beta}^\gamma(q) \quad\textrm{for all }C\in\mathcal S(\gamma,k).$$
Here, $\mathcal S(\alpha,k)$ denotes those $k$-representations belonging to the Segre class $\alpha$. In fact, we may sum over the representations in any two of the three classes and still obtain a polynomial. We illustrate these ideas with an example for the Jordan quiver.

Finally we prove that Hall polynomials exist for all affine quivers with respect to the Bongartz-Dudek decomposition classes \cite{BD}, following the strategy outlined above. We remark that these decomposition classes are also used by Caldero and Reineke in proving the existence of universal polynomials for quiver Grassmannians \cite{CR}. It follows that our result is a refinement of theirs.

In the final section we provide a possible definition of what it should mean for Hall polynomials to exist for a wild quiver. This will necessarily be with respect to some combinatorial partition, and we mention some properties this partition should satisfy.

\section{Representations of quivers}

Let $kQ$ be the path algebra of a connected quiver $Q$ over a field $k$ (see, for example, \cite{ARS}). The category $\modcat kQ$ of finite dimensional $kQ$-modules is an hereditary, abelian, Krull-Schmidt category. It is equivalent to the category of $k$-representations of $Q$, where a $k$-representation is given by a finite dimensional vector space $M_i$ for each vertex $i\in Q_0$ together with a linear map $M_a\colon M_{t(a)}\to M_{h(a)}$ for each arrow $a\in Q_1$. The dimension vector of $M$ is $\dimv M:=(\dim M_i)_i\in\mathbb Z Q_0$ and the Euler form on $\modcat kQ$ satisfies
$$\langle M,N\rangle:=\dim_k\Hom_{kQ}(M,N)-\dim_k\Ext^1_{kQ}(M,N)=[M,N]-[M,N]^1.$$
This depends only on the dimension vectors of $M$ and $N$; indeed
$$\langle\underline d,\underline e\rangle = \sum_{i\in Q_0}d_ie_i-\sum_{a\in Q_1}d_{t(a)}e_{h(a)} \quad\textrm{for all }\underline d,\underline e\in\mathbb ZQ_0.$$
\medskip

A split torsion pair on the set $\indcat kQ$ of indecomposable modules is a decomposition $\indcat kQ=\mathcal F\cup\mathcal T$ such that
$$\Hom(\mathcal T,\mathcal F)=0=\Ext^1(\mathcal F,\mathcal T).$$
The additive subcategories $\mathrm{add}\,\mathcal F$ and $\mathrm{add}\,\mathcal T$ are called the torsion-free and torsion classes respectively. Given a split torsion pair, there exists, for each module $X$, a unique submodule $X_t\leq X$ such that $X_t\in\mathrm{add}\,\mathcal T$ and $X_f:=X/X_t\in\mathrm{add}\,\mathcal F$. Hence $X\cong X_f\oplus X_t$.
\medskip

Let $Q$ be a Dynkin quiver. The dimension vector map induces a bijection between the set of isomorphism classes of indecomposable representations and the set of positive roots $\Phi_+$ of the corresponding semisimple complex Lie algebra \cite{Gabriel}. In particular, this description is independent of the field $k$. Moreover, each indecomposable representation $X$ is a brick ($\End(X)\cong k$) and rigid ($\Ext^1(X,X)=0$), thus exceptional.
\medskip

Let $Q$ be the Jordan quiver $Q$, having one vertex and one loop. Then $kQ=k[t]$ is a principal ideal domain, so finite dimensional modules are described by their elementary divisors. In particular we can associate to a finite dimensional module $M$ the data $\{(\lambda_1,p_1),\ldots,(\lambda_r,p_r)\}$ consisting of partitions $\lambda_i$ and distinct monic irreducible polynomials $p_i\in k[t]$ such that
$$M \cong \bigoplus_{i=1}^r M(\lambda_i,p_i),$$
where, for a partition $\lambda=(1^{l_1}\cdots n^{l_n})$ and monic irreducible polynomial $p$, we write
$$M(\lambda,p)=\bigoplus_r\big(k[t]/(p^r)\big)^{l_r}.$$
\medskip

Let $Q$ be an oriented cycle. Then the category $\modcat^0kQ$ of nilpotent modules is uniserial, with simple modules parameterised by the vertices of $Q$ \cite{Ringel3}. Each indecomposable is determined by its simple top and Loewy length, so the set of all isomorphism classes is in bijection with support-finite functions $(Q_0\times\mathbb N)\to\mathbb N_0$.
\medskip

Let $Q$ be an extended Dynkin quiver which is not an oriented cycle. The roots of $Q$ are either real or imaginary, $\Phi=\Phi^{\mathrm{re}}\cup\Phi^{\mathrm{im}}$, with each imaginary root a non-zero integer multiple of a positive imaginary root $\delta$ \cite{Kac2}. In studying indecomposable representations, Dlab and Ringel \cite{DR} showed the importance of the defect map
$$\partial\colon\mathbb Z Q_0\to\mathbb Z, \quad \underline e\mapsto\langle\delta,\underline e\rangle.$$
By definition, this map is additive on short exact sequences.

We call an indecomposable $kQ$-module $M$ preprojective if $\partial(M)<0$, preinjective if $\partial(M)>0$ and regular if $\partial(M)=0$. This yields a decomposition of $\indcat kQ$ into a `split torsion triple', $\indcat kQ=\mathcal P\cup\mathcal R\cup\mathcal I$, where $\mathcal P$ is the set of indecomposable preprojective modules, $\mathcal I$ the set of indecomposable preinjective modules and $\mathcal R$ the set of indecomposable regular modules. In particular,
\begin{align*}
\Hom(\mathcal I,\mathcal R) &= 0, & \Hom(\mathcal I,\mathcal P) &= 0, & \Hom(\mathcal R,\mathcal P) &= 0,\\
\Ext^1(\mathcal R,\mathcal I) &= 0, & \Ext^1(\mathcal P,\mathcal I) &= 0, & \Ext^1(\mathcal P,\mathcal R) &= 0.
\end{align*}
The indecomposable preprojective and preinjective modules are exceptional. In particular, these indecomposables are determined up to isomorphism by their dimension vectors. Moreover, there is a partial order $\succeq$ on the set of indecomposable preinjective modules such that
$$\Hom(M,N)\neq0 \quad\textrm{implies}\quad M\succeq N.$$
The minimal elements are the simple injective modules.

On the other hand, the category of regular modules is an abelian, exact subcategory which decomposes into a direct sum of serial categories, or tubes. Each tube has a finite number of quasi-simples, forming a single orbit under the Auslander-Reiten translate (say of size $p$), and thus is equivalent to the category of nilpotent representations of a cyclic quiver (with $p$ vertices) over a finite field extension of $k$ (given by the endomorphism ring of any quasi-simple in the tube). Moreover, there are at most three non-homogeneous tubes (those having period $p\geq2$), and the corresponding quasi-simples are exceptional.

Since each tube is a serial category, we see that the isomorphism class of a module without homogeneous regular summands can be described combinatorially.
\medskip

A pair $(A,B)$ of modules is called an orthogonal exceptional pair provided that $A$ and $B$ are exceptional modules such that
$$\Hom(A,B)=\Hom(B,A)=\Ext^1(B,A)=0.$$
We denote by $\mathcal F(A,B)$ the full subcategory of objects having a filtration with factors $A$ and $B$. Then $\mathcal F(A,B)$ is an exact, hereditary, abelian subcategory equivalent to the category of modules over the finite dimensional, hereditary $k$-algebra
$$\begin{pmatrix} k&k^d\\0&k\end{pmatrix} = k(\cdot\xrightarrow{d}\cdot),$$
where $d=\dim\Ext^1(A,B)$ and the quiver above has $d$ arrows from left to right.

\begin{Thm}[Schofield, \cite{Ringel5}]\label{Schofield}
If $M$ is exceptional but not simple, then $M\in\mathcal F(A,B)$ for some orthogonal exceptional pair $(A,B)$, and $M$ is not a simple object in $\mathcal F(A,B)$. In fact, there are precisely $s(M)-1$ such pairs, where $s(M)$ is the size of the support of $\dimv M$.
\end{Thm}

\begin{Lem}\label{Lem_Schofield}
Let $(A,B)$ be an orthogonal exceptional pair of $kQ$-modules and set $d:=[A,B]^1$. If $Q$ is Dynkin, then $d\leq1$, and if $Q$ is extended Dynkin, then $d\leq 2$. Moreover, if $Q$ is extended Dynkin and $d=2$, then $\dimv(A\oplus B)=\delta$ and $\partial(A)=1$.
\end{Lem}

\begin{proof}
If $d\geq3$, then the algebra $\begin{pmatrix} k&k^d\\0&k\end{pmatrix}$ is wild, whereas if $d=2$, then this algebra is tame. The first result follows.

Now let $Q$ be extended Dynkin and $d=2$. There are two indecomposable modules $R\not\cong R'$ which are both extensions of $A$ by $B$. Using the description of $\indcat kQ$ above, we see that $\dimv(A\oplus B)=\dimv R=r\delta$ for some $r\geq1$. Now
$$1=\langle A\oplus B, A\rangle=\langle r\delta,\dimv A\rangle=r\partial(A).$$
Hence $r=\partial(A)=1$.
\end{proof}

Let $Q$ be extended Dynkin. We fix a preprojective module $P$ and a preinjective module $I$ such that $\dimv P+\dimv I=\delta$ and $\partial(I)=1$. Note that $\partial(P)=-1$ and that $P$ and $I$ are necessarily indecomposable. It follows that $(I,P)$ is an orthogonal exceptional pair such that $d=[I,P]^1=2$. Thus there is an embedding of the module category $\modcat K$ of the Kronecker algebra
$$K:=\begin{pmatrix}k&k^2\\0&k\end{pmatrix}=k(\cdot\rightrightarrows\cdot)$$
into $\modcat kQ$ which sends the simple projective to $P$ and the simple injective to $I$.

Under this embedding we can identify the tubes of $\modcat kQ$ with those of $\modcat K$, and it is well-known that the tubes of $\modcat K$ are parameterised by the scheme-theoretic points of the projective line. Moreover, we may assume that the non-homogeneous tubes correspond to some subset of $\{0,1,\infty\}$, and if $R(x)\in\modcat kQ$ is the quasi-simple in the homogeneous tube labelled by $x\in\mathbb P^1_k$, then $\End(R(x))\cong\kappa(x)$ is given by the residue field at $x$ and $\dimv R(x)=m\delta$ where $m=\deg x=[\kappa(x):k]$.

In fact, since the indecomposable preprojective and preinjective modules are uniquely determined by their dimension vectors, it is possible to take a closed subscheme $\mathbb H_{\mathbb Z}\subset\mathbb P^1_{\mathbb Z}$, defined over the integers, such that for any field $k$, the scheme $\mathbb H_k$ parameterises the homogeneous tubes in $\modcat kQ$.

It follows that homogeneous regular modules can be described by pairs consisting of a partition together with a point of the scheme $\mathbb H_k$.

\section{Hall algebras of quiver representations}

Now let $k$ be a finite field. In \cite{Ringel1} Ringel introduced the Ringel-Hall algebra $\HA=\HA(\modcat kQ)$, a free abelian group with basis the set of isomorphism classes of finite dimensional $kQ$-modules and multiplication
$$[M][N]:=\sum_{[X]}F_{MN}^X[X].$$
The structure constants $F_{MN}^X$ are called Hall numbers and are given by
$$F_{MN}^X:=|\{U\leq X:U\cong N, X/U\cong M\}|=\frac{P_{MN}^X}{a_Ma_N},$$
where $a_M:=|\Aut(M)|$ and $P_{MN}^X:=|\{(f,g):0\to N\xrightarrow{f} X \xrightarrow{g} M\to 0\textrm{ exact}\}|$. We also have Riedtmann's Formula \cite{Riedtmann}
$$F_{MN}^X=\frac{|\Ext^1(M,N)_X|}{|\Hom(M,N)|}\cdot\frac{a_X}{a_Ma_N},$$
where $\Ext^1(M,N)_X$ is the set of classes of extensions of $M$ by $N$ with middle term isomorphic to $X$.

Dually, there is a natural comultiplication
$$\Delta\colon\HA\to\HA\otimes\HA, \qquad \Delta([X])=\sum_{[M],[N]}\frac{P_{MN}^X}{a_X}[M]\otimes[N],$$
and $\HA$ is both an associative algebra with unit $[0]$, and a coassociative coalgebra with counit $\epsilon([M])=\delta_{M0}$. Both these statements follow from the identity
$$\sum_{[X]}F_{AB}^XF_{XC}^M = \sum_{[X]}F_{AX}^MF_{BC}^X,$$
a consequence of the pull-back/push-out constructions.

A milestone in the theory of Ringel-Hall algebras was the proof by Green that the Ringel-Hall algebra is a twisted bialgebra (\cite{Green, Ringel4}). That is, we define a new multiplication on the tensor product $\HA\otimes\HA$ via
$$([A]\otimes[B])\cdot([C]\otimes[D]):=q^{\langle A,D\rangle}[A][C]\otimes[B][D].$$
Then the multiplication and comultiplication are compatible with respect to this new multiplication on the tensor product; that is,
$$\Delta([M][N])=\Delta([M])\cdot\Delta([N]).$$
The proof reduces to Green's Formula
$$\sum_{[E]}F_{MN}^EF_{XY}^E/a_E = \sum_{[A],[B],[C],[D]}q^{-\langle A,D\rangle}F_{AB}^MF_{CD}^NF_{AC}^XF_{BD}^Y\frac{a_Aa_Ba_Ca_D}{a_Ma_Na_Xa_Y}.$$
There is also a positive-definite pairing on the Ringel-Hall algebra, and (after twisting the multiplication by the Euler form) the Ringel-Hall algebra is naturally isomorphic to the specialisation at $v=\sqrt q$ of the quantised enveloping algebra of the positive part of a (generalised) Kac-Moody Lie algebra \cite{SVdB}. From this one can deduce part of Kac's Theorem that the set of dimension vectors of indecomposable modules coincides with the set of positive roots $\Phi_+$ of the quiver \cite{DX, Kac1}.
\medskip

Suppose that $M$ is exceptional. Then $\Aut(M^m)\cong\mathrm{GL}_m(k)$ and Riedtmann's Formula implies that
$$F_{MM^m}^{M^{m+1}} = q^{-m\langle M,M\rangle}\frac{a_{M^{m+1}}}{a_Ma_{M^m}} = \frac{q^{m+1}-1}{q-1}.$$
Alternatively we can use the quiver Grassmiannian $\mathrm{Gr}\binom{M^{m+1}}{M}\cong\mathbb P^m_k$, which parameterises embeddings $M\hookrightarrow M^{m+1}$.
\medskip

If $(\mathcal F,\mathcal T)$ is a split torsion pair and $X\cong X_f\oplus X_t$ a module with $X_f\in\mathrm{add}\,\mathcal F$ and $X_t\in\mathrm{add}\,\mathcal T$, then
$$F_{X_fX_t}^E = \delta_{EX} \quad\textrm{and}\quad a_X = a_{X_f}a_{X_t}q^{\langle X_f,X_t\rangle}.$$
These observations will be used repeatedly.

\section{Representation-finite hereditary algebras}

Let $Q$ be a Dynkin quiver and denote by $\Phi_+$ the set of positive roots of the corresponding finite dimensional, semisimple complex Lie algebra. Recall that for any field $k$ the isomorphism classes of $kQ$-modules are in bijection with functions $\alpha\colon\Phi_+\to\mathbb N_0$.

In this section, we provide a new proof of the existence of Hall polynomials for Dynkin quivers. The proof relies on combining Green's Formula with split torsion pairs. We begin with an easy lemma.
\begin{Lem}
For each function $\alpha\colon\Phi_+\to\mathbb N_0$ there exists a monic polynomial $a_\alpha\in\mathbb Z[T]$ such that, for any finite field $k$ with $|k|=q$ and any $k$-representation $A$ with $[A]=\alpha$,
$$a_\alpha(q)=a_A=|\Aut(A)|.$$
\end{Lem}
\begin{proof}
Each indecomposable is exceptional and the dimension of homomorphisms between indecomposables is given via the Auslander-Reiten quiver and the mesh relations, hence is independent of the field $k$.
\end{proof}
Given a split torsion pair $(\mathcal F,\mathcal T)$ and a module $X\cong X_f\oplus X_t$, we can simplify the left hand side of Green's Formula as follows:
$$\sum_EF_{MN}^EF_{X_fX_t}^E/a_E=F_{MN}^X/a_X=F_{MN}^X/a_{X_f}a_{X_t}q^{\langle X_f,X_t\rangle}.$$
(Note that the sum is actually over isomorphism classes of representations, though we shall often use this more convenient notation.) If we now consider the right hand side of Green's Formula, we see that all Hall numbers involve middle terms with dimension vector strictly smaller than $\dimv X$ (provided that $M$, $N$, $X_f$ and $X_t$ are all non-zero, of course):
$$\sum_{A,B,C,D}q^{-\langle A,D\rangle}F_{AB}^MF_{CD}^NF_{AC}^{X_f}F_{BD}^{X_t}\frac{a_Aa_Ba_Ca_D}{a_Ma_Na_{X_f}a_{X_t}}.$$
Thus
$$F_{MN}^X = \sum_{A,B,C,D}q^{\langle X_f,X_t\rangle-\langle A,D\rangle}F_{AB}^MF_{CD}^NF_{AC}^{X_f}F_{BD}^{X_t}\frac{a_Aa_Ba_Ca_D}{a_Ma_N}.$$
\begin{Thm}[Ringel \cite{Ringel2}]\label{Ringel}
Hall polynomials exist for $Q$; that is, given $\mu$, $\nu$ and $\xi$, there exists an integer polynomial $F_{\mu\nu}^\xi\in\mathbb Z[T]$ such that, for any finite field $k$ with $|k|=q$ and any $k$-representations $M$, $N$ and $X$ with $[M]=\mu$, $[N]=\nu$ and $[X]=\xi$, we have
$$F_{\mu\nu}^{\xi}(q) = F_{MN}^X.$$
\end{Thm}

\begin{proof}
We wish to show that $F_{MN}^X$ is given by some universal integer polynomial. Let $X=I_1^{m_1}\oplus\cdots I_r^{m_r}$ be a decomposition of $X$ into pairwise non-isomorphic indecomposable modules $I_i$. Then (up to reordering) there exists a split torsion pair such that $X_t=I_r^{m_r}$. We now use the formula above together with induction on dimension vector to deduce that $F_{MN}^X$ is of the form (polynomial)/(monic polynomial). Since this must take integer values at all prime powers, we have that $F_{MN}^X$ is given by some universal polynomial $F_{MN}^X(T)\in\mathbb Z[T]$. It is thus enough to consider the case when $X=I^m$ is isotypic.

Suppose $X=I^{m+1}$ with $I$ indecomposable. Since $I$ is exceptional, we can simplify the left hand side of Green's Formula as
$$\sum_EF_{MN}^EF_{II^m}^E/a_E = F_{MN}^{I^{m+1}}/q^ma_Ia_{I^m},$$
whence
$$F_{MN}^{I^{m+1}} = \sum_{A,B,C,D}q^{m-\langle A,D\rangle}F_{AB}^MF_{CD}^NF_{AC}^IF_{BD}^{I^m}\frac{a_Aa_Ba_Ca_D}{a_Ma_N}.$$
By induction on dimension vector, $F_{MN}^{I^{m+1}}$ is given by a universal integer polynomial. We thus reduce to the case when $X=I$ is indecomposable.
Gabriel's Theorem tells us that the indecomposable modules are determined up to isomorphism by their dimension vectors, so if $\dimv E=\dimv X$ and $E\not\cong X$, then $E$ must be decomposable. Set $R:=\mathrm{rad}\,X$ and $T:=X/R$, so that $F_{TR}^X=1$. Then
$$F_{MN}^X/a_X = \sum_EF_{MN}^EF_{TR}^E/a_E - \sum_{E\,\mathrm{decomp}}F_{MN}^EF_{TR}^E/a_E.$$
The second sum on the right hand side is over decomposable modules, hence is of the form (polynomial)/(monic polynomial), as is the first sum by Green's Formula and induction on dimension vector. Thus $F_{MN}^X$ is given by a universal integer polynomial.

Alternatively, Theorem \ref{Schofield} gives an exact sequence $0\to B\to X\to A\to 0$ for some orthogonal exceptional pair $(A,B)$ with $[A,B]^1=1$. Either way, we reduce to the case when $X$ is simple, where the result is trivial.
\end{proof}

\section{Cyclic quivers}

Let $Q$ be an oriented cycle and $k$ an arbitrary field. Recall that the set of isomorphism classes of nilpotent modules is in bijection with support-finite functions $(Q_0\times\mathbb N)\to\mathbb N_0$.
It is known from \cite{Ringel3} that Hall polynomials exist in this context. We provide a new proof of this fact, using Green's Formula and induction on partitions. Given a module $X$ denote by $\lambda(X)$ the partition formed by taking the Loewy lengths of its indecomposable summands. We order partitions via the reverse lexicogarphic ordering; that is, if $\lambda=(1^{l_1}2^{l_2}\cdots n^{l_n})$ and $\mu=(1^{m_1}2^{m_2}\cdots n^{m_n})$ are written in exponential form, then
$$\lambda<\mu \quad\textrm{if there exists $i$ such that $l_i>m_i$ and $l_j=m_j$ for all $j>i$.}$$
The cup product $\lambda\cup\mu$ is the partition $(1^{l_1+m_1}\cdots n^{l_n+m_n})$.
As in the representation-finite case, we have the following lemma.
\begin{Lem}
The dimension $\dim\Hom(A,B)$ depends only on $\alpha=[A]$ and $\beta=[B]$ and not on the choice of field $k$. Moreover, there exists a monic polynomial $a_\alpha\in\mathbb Z[T]$ such that, for any finite field $k$ with $|k|=q$ and any $k$-representation $A$ with $[A]=\alpha$, we have $a_\alpha(q)=a_A$.
\end{Lem}
\begin{proof}
Since indecomposables are uniserial, it is clear that the dimension of the space of homomorphisms between indecomposables is independent of the field. The lemma follows easily.
\end{proof}
\begin{Lem}[\cite{Hubery}]\label{cyclic-ext}
If $X$ is an extension of $M$ by $N$, then $\lambda(X)\leq\lambda(M)\cup\lambda(N)$ with equality if and only if $X\cong M\oplus N$.
\end{Lem}
\begin{Thm}\label{cyclic}
Hall polynomials exist for nilpotent representations of cyclic quivers.
\end{Thm}

\begin{proof}
Let $X=X_1\oplus X_2$ be decomposable. Then
$$F_{MN}^XF_{X_1X_2}^X/a_X = \sum_EF_{MN}^EF_{X_1X_2}^E/a_E - \sum_{\lambda(E)<\lambda(X)}F_{MN}^EF_{X_1X_2}^E/a_E.$$
By induction on $\lambda(X)$ we know that the second sum on the right hand side is of the form (polynomial)/(monic polynomial), as is the first sum by Green's Fromula and induction on dimension vector. Moreover, by Riedtmann's Formula, we know that
$$F_{X_1X_2}^X/a_X=1/q^{[X_1,X_2]}a_{X_1}a_{X_2},$$
the reciprocal of a monic polynomial. Hence $F_{MN}^X$ is given by a universal integer polynomial.
We are reduced to proving the formula when $X$ is indecomposable. Since the category is uniserial, this implies that both $M$ and $N$ are indecomposable, in which case the Hall number is either $1$ or $0$ and is independent of the field.
\end{proof}
\textbf{Remark.}
The same induction can be used to show that
$$2\deg F_{MN}^X\leq \deg a_X-\deg(a_Ma_N).$$
For, if $X\cong X_1\oplus X_2$, then induction gives
$$2\deg F_{MN}^X\leq \deg(a_{X_1}a_{X_2})+2[X_1,X_2]-\deg(a_Ma_N),$$
and obviously
$$\deg a_X = \deg(a_{X_1}a_{X_2})+[X_1,X_2]+[X_2,X_1].$$
Hence taking such a decomposition with $[X_1,X_2]\leq[X_2,X_1]$ gives the result.

In fact, for the Jordan quiver, we always have equality and the leading coefficient is given by the corresponding Littlewood-Richardson coefficient \cite{Macdonald}.
\medskip

Consider now the Jordan quiver $Q$, so that $kQ=k[t]$. Recall that a finite dimensional module $M$ is determined by the data $\{(\lambda_1,p_1),\ldots,(\lambda_r,p_r)\}$ consisting of partitions $\lambda_i$ and distinct monic irreducible polynomials $p_i\in k[t]$ such that
$$M \cong \bigoplus_{i=1}^r M(\lambda_i,p_i).$$
Clearly the primes $p_i$ depend on the field, but we can partition the set of isomorphism classes by considering just their degrees. This is called the Segre decomposition. More precisely, a Segre symbol is a multiset $\sigma=\{(\lambda_1,d_1),\ldots,(\lambda_r,d_r)\}$ of pairs $(\lambda,d)$ consisting of a partition $\lambda$ and a positive integer $d$. The corresponding Segre class $\mathcal S(\sigma,k)$ consists of those isomorphism classes of modules of type $\{(\lambda_1,p_1),\ldots,(\lambda_r,p_r)\}$, where the $p_i\in k[t]$ are distinct monic irreducible polynomials with $\deg p_i=d_i$.
\begin{Thm}[\cite{Arnold,BD,Gibson}]
Let $k$ be an algebraically closed field. Then the Segre classes stratify the variety $\End(k^m)$ into smooth, irreducible, $\mathrm{GL}_m(k)$-stable subvarieties, each admitting a smooth, rational, geometric quotient. Moreover, the stabilisers of any two matrices in the same Segre class are conjugate inside $\mathrm{GL}_m(k)$.
\end{Thm}
\begin{Lem}
Given a Segre symbol $\sigma$, there exists a monic polynomial $a_\sigma\in\mathbb Z[T]$ and a polynomial $n_\sigma\in\mathbb Q[T]$ such that, for any finite field $k$ with $|k|=q$,
$$a_\sigma(q)=a_M \quad\textrm{for any }M\in\mathcal S(\sigma,k), \textrm{ and}\quad n_\sigma(q)=|\mathcal S(\sigma,k)|.$$
\end{Lem}
\begin{proof}
If $\sigma=\{(\lambda_1,d_1),\ldots,(\lambda_r,d_r)\}$, then $a_\sigma(T)=\prod_ia_{\lambda_i}(T^{d_i})$ where, for a partition $\lambda=(1^{l_1}\cdots n^{l_n})$,
$$a_\lambda(T):=T^{\sum_{i,j}\min\{i,j\}l_il_j}\prod_i(1-T^{-1})\cdots(1-T^{-l_i}).$$
This polynomial occurs in \cite{Macdonald} as the size of the automorphism group of the module $M(\lambda,t)$.

We write $\sigma(d)$ to be the Segre symbol formed by those pairs $(\lambda_i,d_i)$ in $\sigma$ with $d_i=d$. To obtain the formula for $n_\sigma$, let us first suppose that $\sigma=\sigma(d)$ for some $d$. Then
$$n_\sigma=\phi_d(\phi_d-1)\cdots(\phi_d-r+1)/z_\sigma,$$
where $\phi_d$ is the number of monic irreducible polynomials of degree $d$ and $z_\sigma$ is the size of the stabiliser for the natural action of the symmetric group $\mathfrak S_r$ on $(\lambda_1,\ldots,\lambda_r)$ given by place permutation.

In general we can write $\sigma=\bigcup_d\sigma(d)$, and $n_\sigma=\prod_dn_{\sigma(d)}$.
\end{proof}
\begin{Thm}\label{gen_cyclic}
Given three Segre symbols $\rho$, $\sigma$ and $\tau$, there exists an integer polynomial $F_{\rho\sigma}^\tau\in\mathbb Z[T]$ such that, for any finite field $k$ with $|k|=q$,
$$F_{\rho\sigma}^\tau(q)=\sum_{\substack{R\in\mathcal S(\rho,k)\\S\in\mathcal S(\sigma,k)}}F_{RS}^T \quad\textrm{for all }T\in\mathcal S(\tau,k).$$
Moreover, we have the identities
\begin{align*}
n_\tau(q) F_{\rho\sigma}^\tau(q) &= n_\rho(q)\sum_{\substack{S\in\mathcal S(\sigma,k)\\T\in\mathcal S(\tau,k)}}F_{RS}^T \quad\textrm{for all }R\in\mathcal S(\rho,k)\\
&= n_\sigma(q)\sum_{\substack{R\in\mathcal S(\rho,k)\\T\in\mathcal S(\tau,k)}}F_{RS}^T \quad\textrm{for all }S\in\mathcal S(\sigma,k).
\end{align*}
It follows that both $\frac{n_\tau}{n_\rho}F_{\rho\sigma}^\tau$ and $\frac{n_\tau}{n_\sigma}F_{\rho\sigma}^\tau$ are polynomials (but with rational coefficients).
\end{Thm}
\begin{proof}
The proof is similar to when we just considered nilpotent modules (i.e. the single irreducible polynomial $p(t)=t$). We begin by noting that there are no homomorphisms between modules corresponding to distinct irreducible polynomials. In particular we can decompose $\rho=\bigcup_d\rho(d)$, and given $R\in\mathcal S(\rho,k)$, there is a unique decomposition $R=\bigoplus_dR(d)$ such that $R(d)\in\mathcal S(\rho(d),k)$. We deduce that
$$F_{RS}^T = \prod_dF_{R(d)S(d)}^{T(d)}.$$
In particular, we can reduce to the case when all degrees which occur in $\rho$, $\sigma$ and $\tau$ equal some fixed integer $d$.

Let us fix
$$T:=M(\nu_1,p_1)\oplus\cdots\oplus M(\nu_m,p_m)\in\mathcal S(\tau,k).$$
By adding in copies of the zero partition, we may assume that $\rho=(\lambda_1,\ldots,\lambda_m)$ and $\sigma=(\mu_1,\ldots,\mu_m)$ for the same $m$ (we have simplified the notation by omitting the number $d$). Since there are no homomorphisms between modules corresponding to distinct irreducible polynomials, if $F_{RS}^T$ is non-zero for some $R\in\mathcal S(\rho,k)$ and $S\in\mathcal S(\sigma,k)$, then there exist permutations $r,s\in\mathfrak S_m$ such that
$$R=\bigoplus_iM(\lambda_{r(i)},p_i), \quad S=\bigoplus_iM(\mu_{s(i)},p_i) \quad\textrm{and}\quad F_{RS}^T=\prod_i F_{\lambda_{r(i)}\mu_{s(i)}}^{\nu_i}(q^d),$$
where $F_{\lambda\mu}^\nu$ is the classical Hall polynomial.

It follows that
$$\sum_{R\in\mathcal S(\rho), S\in\mathcal S(\sigma)}F_{RS}^T = \sum_{r,s}\prod_i F_{\lambda_{r(i)}\mu_{s(i)}}^{\nu_i}(q^d),$$
where the sum is taken over all permutations $r$ and $s$ yielding non-isomorphic modules; that is, $r$ runs through the cosets in $\mathfrak S_m$ with respect to the stabiliser of $(\lambda_1,\ldots,\lambda_m)$, and similarly for $s$. It is now clear that this number is described by a universal polynomial over the integers, which we denote by $F_{\rho\sigma}^\tau$.
Suppose instead that we fix
$$R=M(\lambda_1,p_1)\oplus\cdots\oplus M(\lambda_m,p_m)\in\mathcal S(\rho),$$
where $\lambda_1,\ldots,\lambda_m$ are all non-zero.

Assume first that $\tau=\{\nu_1,\ldots,\nu_m\}$ consists of precisely $m$ partitions. By adding in copies of the zero partition, we may further assume that $\sigma=\{\mu_1,\ldots,\mu_m\}$ also consists of $m$ partitions. It follows as before that
$$\sum_{S\in\mathcal S(\sigma),T\in\mathcal S(\tau)}F_{RS}^T = \sum_{s,t}\prod_iF_{\lambda_i\mu_{s(i)}}^{\nu_{t(i)}}(q^d).$$
In general, $\tau$ will consist of $m+n$ partitions, and if $F_{RS}^T\neq0$, then we can write $T=T'\oplus X$ and $S=S'\oplus X$ such that $T'$ contains all summands of $T$ corresponding to the polynomials $p_i$ occurring in $R$. We observe that $F_{RS}^T=F_{RS'}^{T'}$. It follows that
$$\sum_{\substack{S\in\mathcal S(\sigma)\\T\in\mathcal S(\tau)}}F_{RS}^T = \sum_{\substack{\sigma=\sigma'\cup\xi\\\tau=\tau'\cup\xi}}N_\xi\sum_{\substack{S'\in\mathcal S(\sigma')\\T'\in\mathcal S(\tau')}}F_{RS'}^{T'}.$$
The number $N_\xi$ equals the number of isomorphism classes of $X=\bigoplus_iM(\xi_i,x_i)\in\mathcal S(\xi)$ such that the polynomials $p_1,\ldots,p_m,x_1,\ldots,x_n$ are pairwise distinct. Thus
$$N_\xi=(\phi_d-m)\cdots(\phi_d-m-n+1)/z_\xi$$
and hence the number $\sum_{S,T}F_{RS}^T$ is again given by a universal polynomial. It follows immediately that
$$n_\tau\sum_{R,S}F_{RS}^T=\sum_{R,S,T}F_{RS}^T=n_\rho\sum_{S,T}F_{RS}^T.$$
An analogous argument works for the sum $\sum_{R,T}F_{RS}^T$.
\end{proof}

\section{An example}

Theorem \ref{gen_cyclic} gives a good generalisation of Hall polynomials in the case of arbitrary $k[t]$-modules. In this section we illustrate the proof of Theorem \ref{gen_cyclic} and in so doing, show that it is not possible to fix two modules and still obtain a universal polynomial.

We fix the degree $d=1$ and the Segre symbols
$$\rho:=\{(1,1),(1,1,1),(2,1)\}, \quad \sigma:=\{(1),(1)\} \quad\textrm{and}\quad \tau:=\{(1,1,1),(2,1,1),(2,1)\}.$$
Note that
$$n_\sigma=\frac{1}{2}q(q-1) \quad\textrm{and}\quad n_\rho=n_\tau=q(q-1)(q-2).$$
We are interested in the numbers $F_{RS}^T$ for the modules
\begin{align*}
R &:= M((1,1),x)\oplus M((1,1,1),y)\oplus M((2,1),z)\\
S &:= M((1),x')\oplus M((1),y')\\
T &:= M((1,1,1),x'')\oplus M((2,1,1),y'')\oplus M((2,1),z''),
\end{align*}
where $x$, $y$ and $z$ are distinct elements of $k$, as are $x'$ and $y'$, and $x''$, $y''$ and $z''$.

We begin by computing the possible Hall polynomials that can appear. We have
$$F_{(1,1)(1)}^{(1,1,1)} = q^2+q+1 \quad\textrm{and}\quad F_{(2,1)(1)}^{(2,1,1)} = q(q+1).$$
Both of these sequences are split, so we can apply Riedtmann's Formula. Also
$$F_{(1,1)(1)}^{(2,1)} = 1 \quad\textrm{and}\quad F_{(1,1,1)(1)}^{(2,1,1)} = 1.$$
These are clear, since in both cases we just have the top and radical of the middle term.

We deduce that
$$F_{RS}^T = \begin{cases}
q^2+q+1 &\textrm{if $(x,y,z)=(x'',y'',z'')$ and $\{x',y'\}=\{x,y\}$;}\\
q^2+q   &\textrm{if $(x,y,z)=(z'',x'',y'')$ and $\{x',y'\}=\{x,z\}$;}\\
0       &\textrm{otherwise.}\end{cases}$$

Hence
\begin{align*}
\sum_R F_{RS}^T &= \begin{cases}
q^2+q+1 &\textrm{if $\{x',y'\}=\{x'',y''\}$;}\\
q^2+q   &\textrm{if $\{x',y'\}=\{y'',z''\}$;}\\
0       &\textrm{otherwise,}\end{cases}\\
{}
\sum_S F_{RS}^T &= \begin{cases}
q^2+q+1 &\textrm{if $(x,y,z)=(x'',y'',z'')$;}\\
q^2+q   &\textrm{if $(x,y,z)=(z'',x'',y'')$;}\\
0       &\textrm{otherwise,}\end{cases}\\
{}
\sum_T F_{RS}^T &= \begin{cases}
q^2+q+1 &\textrm{if $\{x',y'\}=\{x,y\}$;}\\
q^2+q   &\textrm{if $\{x',y'\}=\{x,z\}$;}\\
0       &\textrm{otherwise.}\end{cases}
\end{align*}

Finally
\begin{align*}
\sum_{R,S} F_{RS}^T &= 2q^2+2q+1,\\
\sum_{S,T} F_{RS}^T &= 2q^2+2q+1,\\
\sum_{R,T} F_{RS}^T &= 2(2q^2+2q+1)(q-2).
\end{align*}
Thus we only get universal polynomials if we sum over two of the Segre classes.

\section{Tame hereditary algebras}

Let $Q$ be an extended Dynkin quiver which is not an oriented cycle and $k$ a finite field. Recall that we have the `split torsion triple' $\indcat kQ=\mathcal P\cup\mathcal R\cup\mathcal I$ given by the indecomposable preprojective, regular and preinjective modules. Moreover, the category of regular modules decomposes into a direct sum of tubes indexed by the projective line in such a way that each regular simple module $R$ in the tube labelled by $x$ satisfies $\End(R)\cong\kappa(x)$ and $\dimv R=(\deg x)\delta$. Finally, we may also assume that the non-homogeneous tubes are labelled by some subset of $\{0,1,\infty\}$, whereas the homogeneous tubes are labelled by the points of the scheme $\mathbb H_{\mathbb Z}\otimes k$ for some closed integral subscheme $\mathbb H_{\mathbb Z}\subset\mathbb P^1_{\mathbb Z}$.

The indecomposable preprojective and preinjective modules are all exceptional, as are the regular simple modules in the non-homogeneous tubes. Hence the isomorphism class of a module without homogeneous regular summands can be described combinatorially, whereas homogeneous regular modules are determined by pairs consisting of a partition together with a point of the scheme $\mathbb H_k$.
\medskip

We are now in a position to define the partition of Bongartz and Dudek \cite{BD}. A decomposition symbol is a pair $\alpha=(\mu,\sigma)$ such that $\mu$ specifies a module without homogeneous regular summands and $\sigma=\{(\lambda_1,d_1),\ldots,(\lambda_r,d_r)\}$ is a Segre symbol.

Given a decomposition symbol $\alpha=(\mu,\sigma)$ and a field $k$, we define $\mathcal S(\alpha,k)$ to be the set of isomorphism classes of modules of the form $M(\mu,k)\oplus R$, where $M(\mu,k)$ is the $kQ$-module determined by $\mu$ and
$$R=R(\lambda_1,x_1)\oplus\cdots\oplus R(\lambda_r,x_r)$$
for some distinct points $x_1,\ldots,x_r\in\mathbb H_k$ such that $\deg x_i=d_i$.

\begin{Thm}[Bongartz-Dudek \cite{BD}]
If $k$ is algebraically closed, each decomposition class $\mathcal S(\alpha,k)$ determines a smooth, irreducible, $\mathrm{GL}$-invariant subvariety of the corresponding representation variety, which furthermore admits a smooth, rational geometric quotient.
\end{Thm}

Note that it is still open as to whether the closure of a decomposition class is again the union of decomposition classes. Also, unlike in the classical case, the endomorphism algebras of modules in the same decomposition class are not necessarily isomorphic as algebras. We do, however, have the following result.

\begin{Lem}
Given a decomposition symbol $\alpha=(\mu,\sigma)$, there exist universal polynomials $a_\alpha$ and $n_\alpha$ such that, for any finite field $k$ with $|k|=q$,
$$a_\alpha(q)=|\Aut(A)| \quad\textrm{for all }A\in\mathcal S(\alpha,k), \textrm{ and}\quad n_\alpha(q)=|\mathcal S(\alpha,k)|.$$
Moreover, $a_\alpha$ is a monic integer polynomial.
\end{Lem}

We can now state the main result of this paper.
\begin{MainThm}
Hall polynomials exist with respect to the decomposition classes described above; that is, given decomposition classes $\alpha$, $\beta$ and $\gamma$, there exists a rational polynomial $F_{\alpha\beta}^\gamma$ such that, for any finite field $k$ with $|k|=q$,
$$F_{\alpha\beta}^\gamma(q)=\sum_{\substack{A\in\mathcal S(\alpha,k)\\B\in\mathcal S(\beta,k)}}F_{AB}^C \quad\textrm{for all }C\in\mathcal S(\gamma,k)$$
and moreover
\begin{align*}
n_\gamma(q)F_{\alpha\beta}^\gamma(q) &= n_\alpha(q)\sum_{\substack{B\in\mathcal S(\beta,k)\\C\in\mathcal S(\gamma,k)}}F_{AB}^C \quad\textrm{for any }A\in\mathcal S(\alpha,k)\\
&=n_\beta(q)\sum_{\substack{A\in\mathcal S(\alpha,k)\\C\in\mathcal S(\gamma,k)}}F_{AB}^C \quad\textrm{for any }B\in\mathcal S(\beta,k).
\end{align*}
\end{MainThm}
\textbf{Remarks.}
\begin{enumerate}
\item When $\alpha=(\emptyset,\rho)$, $\beta=(\emptyset,\sigma)$ and $\gamma=(\emptyset,\tau)$, so that $\mathcal S(\alpha,k)$, $\mathcal S(\beta,k)$ and $\mathcal S(\gamma,k)$ contain only homogeneous regular modules, the result follows from Theorem \ref{gen_cyclic}.
\item It is not true that each Hall polynomial has integer coefficients. In fact, we can easily construct a counter-example for the Kronecker quiver. Let $\pi_r$ be the decomposition class corresponding to the indecomposable preprojective $P_r$ and consider the three Segre symbols
$$\sigma_1:=\big\{\big((1),1\big),\big((1),1\big)\big\}, \quad \sigma_2:=\big\{\big((2),1\big)\big\} \quad\textrm{and}\quad \sigma_3:=\big\{\big((1),2\big)\big\}.$$
Consider decomposition classes $\alpha_i:=(\emptyset,\sigma_i)$. As mentioned in the introduction, $F_{RP_0}^{P_2}=1$ for all regular modules $R$ containing at most one indecomposable from each tube, and is zero otherwise. Therefore
$$F_{\alpha_1\pi_0}^{\pi_2} = n_{\alpha_1} = q(q+1)/2, \quad F_{\alpha_2\pi_0}^{\pi_2} = n_{\alpha_2} = q+1, \quad F_{\alpha_3\pi_0}^{\pi_2} = n_{\alpha_3} = q(q-1)/2.$$
This does not contradict Proposition 6.1 of \cite{Reineke} since we are not counting the rational points of any scheme. In fact, let $D_r$ be the locally closed subscheme consisting of those representations $X$ such that $\dimv X=2\delta$ and $\dim\mathrm{rad}\End(X)=r$. The quiver Grassmannian $\mathrm{Gr}\binom{P_2}{P_0}\cong\mathbb P^2$ parameterises all submodules of $P_2$ isomorphic to $P_0$, and decomposes into an open subscheme consisting of those submodules whose cokernel lies in $D_0$ (with $F_{\alpha_1\pi_0}^{\pi_2}+F_{\alpha_3\pi_0}^{\pi_2}=q^2$ rational points) and a closed subscheme consisting of those points whose cokernel lies in $D_1$ (with $F_{\alpha_2\pi_0}^{\pi_2}=q+1$ rational points).
\end{enumerate}

Let us call a decomposition class $\alpha$ discrete if $n_\alpha=1$ (hence $\alpha=(\mu,\emptyset)$). For example, if $\alpha$ is discrete, then the Hall polynomial $F_{\alpha\beta}^\gamma$ satisfies
$$F_{\alpha\beta}^\gamma(q) = \sum_{B\in\mathcal S(\beta,k)}F_{AB}^C \quad\textrm{for all }C\in\mathcal S(\gamma,k).$$
We call a module discrete if it is the unique module in a discrete decomposition class, which is if and only if it contains no homogeneous regular summand.

\section{Reductions using Green's Formula}\label{section:Green}

In this section we show how Green's Formula together with split torsion pairs can be used to set up an induction on dimension vector. This reduces the problem of existence of Hall polynomials to the special case of $F_{MN}^X$ with $X$ regular and either homogeneous or non-homogeneous lying in a single tube.

Let $(\mathcal T,\mathcal F)$ be a split torsion pair and suppose that each indecomposable homogeneous regular module is contained in $\mathcal T$. We can decompose any module $A\cong A_f\oplus A_t$ with $A_f\in\mathrm{add}\,\mathcal F$ and $A_t\in\mathrm{add}\,\mathcal T$, and $A_f$ is discrete. Thus every decomposition class $\alpha$ can be written as $\alpha_f\oplus\alpha_t$ with $\alpha_f$ discrete.

A dual result clearly holds if each indecomposable homogeneous regular module is contained in $\mathcal F$.

\begin{Prop}
Suppose that Hall polynomials exist for all modules of dimension vector smaller than $\underline d$ and let $\xi$ be a decomposition class of dimension vector $\underline d$. Let $(\mathcal T,\mathcal F)$ be a split torsion pair such that $\xi_f$ and $\xi_t$ are both non-zero and such that all indecomposable homogeneous regular modules are contained in either $\mathcal F$ or $\mathcal T$. Then the Hall polynomial $F_{\mu\nu}^\xi$ exists and is given by
$$F_{\mu\nu}^\xi = \sum_{\alpha,\beta,\gamma,\delta}q^{\langle\xi_f,\xi_t\rangle-\langle\alpha,\delta\rangle}F_{\alpha\beta}^\mu F_{\gamma\delta}^\nu F_{\alpha\gamma}^{\xi_f} F_{\beta\delta}^{\xi_t} \frac{a_\alpha a_\beta a_\gamma a_\delta}{a_\mu a_\nu}\frac{n_\mu n_\nu}{n_\alpha n_\beta n_\gamma n_\delta}.$$
\end{Prop}

\begin{proof}
The proof involves analysing the influence of the split torsion pair on Green's Formula. We then sum in a suitable way to deduce the existence of Hall polynomials.

Suppose that all indecomposable homogeneous regular modules are contained in $\mathcal T$ and consider Green's Formula. On the left we have
$$\sum_EF_{MN}^EF_{X_fX_t}^E/a_E = F_{MN}^X/a_X = F_{MN}^X/a_{X_f}a_{X_t}q^{\langle X_f,X_t\rangle}$$
whereas the right hand side reads
$$\sum_{A,B,C,D}q^{-\langle A,D\rangle}F_{AB}^MF_{CD}^NF_{AC}^{X_f}F_{BD}^{X_t}\frac{a_Aa_Ba_Ca_D}{a_Ma_Na_{X_f}a_{X_t}}.$$
We thus have the equality
$$F_{MN}^X = \sum_{A,B,C,D}q^{\langle X_f,X_t\rangle-\langle A,D\rangle}F_{AB}^MF_{CD}^NF_{AC}^{X_f}F_{BD}^{X_t}\frac{a_Aa_Ba_Ca_D}{a_Ma_N}.$$
Observe that $C=C_f$ and $B=B_t$.
Writing $A=A_f\oplus A_t$ and $M=M_f\oplus M_t$, associativity of Hall numbers implies
$$F_{AB_t}^M = \sum_LF_{A_fA_t}^LF_{LB_t}^M = \sum_LF_{A_fL}^MF_{A_tB_t}^L = \sum_{L_t}F_{A_fL_t}^MF_{A_tB_t}^{L_t} = \delta_{A_fM_f}F_{A_tB_t}^{M_t}.$$
Similarly $F_{C_fD}^N=\delta_{D_tN_t}F_{C_fD_f}^{N_f}$, so that
$$F_{MN}^X = \sum_{A_t,B_t,C_f,D_f}q^{\langle X_f,X_t\rangle-\langle A,D\rangle} F_{A_tB_t}^{M_t}F_{C_fD_f}^{N_f}F_{AC_f}^{X_f}F_{B_tD}^{X_t}\frac{a_Aa_{B_t}a_{C_f}a_D}{a_Ma_N},$$
where $A=M_f\oplus A_t$ and $D=D_f\oplus N_t$. Therefore
$$\sum_{\substack{M\in\mathcal S(\mu)\\N\in\mathcal S(\nu)}}\!\!\!F_{MN}^X = \sum_{\substack{\alpha_t,\beta_t\\\gamma_f,\delta_f}}q^{\langle\xi_f,\xi_t\rangle-\langle\alpha,\delta\rangle} \big(\!\!\! \sum_{\substack{A_t\in\mathcal S(\alpha_t)\\M_t\in\mathcal S(\mu_t)}}\!\!\!\!\! F_{A_tB_t}^{M_t} \big) F_{C_fD_f}^{N_f} F_{AC_f}^{X_f} \big(\!\!\! \sum_{\substack{B_t\in\mathcal S(\beta_t)\\N_t\in\mathcal S(\nu_t)}}\!\!\!\!\! F_{B_tD}^{X_t} \big) \frac{a_\alpha a_{\beta_t}a_{\gamma_f}a_\delta}{a_\mu a_\nu},$$
where we have written $\alpha=\mu_f\oplus\alpha_t$ and, for $A_t\in\mathcal S(\alpha_t)$, $A=M_f\oplus A_t$. Similarly $\delta=\delta_f\oplus\nu_t$ and, for $N_t\in\mathcal S(\nu_t)$, $D=D_f\oplus N_t$.
By induction, we have the Hall polynomials $F_{\alpha_t\beta_t}^{\mu_t}$, $F_{\beta_t\delta}^{\xi_t}$, $F_{\alpha\gamma_f}^{\xi_f}$ and $F_{\gamma_f\delta_f}^{\nu_f}$ so that, for any $X\in\mathcal S(\xi)$,
\begin{align*}
\sum_{\substack{M\in\mathcal S(\mu)\\N\in\mathcal S(\nu)}}F_{MN}^X &=
\sum_{\substack{\alpha_t,\beta_t\\\gamma_f,\delta_f}}q^{\langle\xi_f,\xi_t\rangle-\langle\alpha,\delta\rangle} \big(F_{\alpha_t\beta_t}^{\mu_t}n_{\mu_t}/n_{\beta_t}\big) F_{\gamma_f\delta_f}^{\nu_f} \big(F_{\alpha\gamma_f}^{\xi_f}/n_\alpha\big) F_{\beta_t\delta}^{\xi_t} \frac{a_\alpha a_{\beta_t}a_{\gamma_f}a_\delta}{a_\mu a_\nu}\\
&= \sum_{\substack{\alpha_t,\beta_t\\\gamma_f,\delta_f}}q^{\langle\xi_f,\xi_t\rangle-\langle\alpha,\delta\rangle} F_{\alpha_t\beta_t}^{\mu_t} F_{\gamma_f\delta_f}^{\nu_f} F_{\alpha\gamma_f}^{\xi_f} F_{\beta_t\delta}^{\xi_t} \frac{a_\alpha a_{\beta_t}a_{\gamma_f}a_\delta}{a_\mu a_\nu} \frac{n_{\mu_t}}{n_\alpha n_{\beta_t}}.
\end{align*}
We now note that the polynomial $F_{\beta\delta}^{\xi_t}$ is non-zero only if $\beta=\beta_t$. In this case, we also have the identity $F_{\alpha\beta}^\mu = \delta_{\alpha_f\mu_f}F_{\alpha_t\beta_t}^{\mu_t}$ exactly as for modules. Similarly $\gamma=\gamma_f$ and $F_{\gamma\delta}^\nu=\delta_{\nu_t\delta_t}F_{\gamma_f\delta_f}^{\nu_f}$. Hence we can simplify the above expression to get
$$F_{\mu\nu}^\xi = 
\sum_{\alpha,\beta,\gamma,\delta}q^{\langle\xi_f,\xi_t\rangle-\langle\alpha,\delta\rangle} F_{\alpha\beta}^\mu F_{\gamma\delta}^\nu F_{\alpha\gamma}^{\xi_f} F_{\beta\delta}^{\xi_t} \frac{a_\alpha a_\beta a_\gamma a_\delta}{a_\mu a_\nu} \frac{n_\mu}{n_\alpha n_\beta}.$$

Similarly
\begin{align*}
\sum_{\substack{M\in\mathcal S(\mu)\\X\in\mathcal S(\xi)}}\!\!\!F_{MN}^X &= \sum_{\substack{\alpha_t,\beta_t\\\gamma_f,\delta_f}}q^{\langle\xi_f,\xi_t\rangle-\langle\alpha,\delta\rangle} \big(\!\!\! \sum_{\substack{A_t\in\mathcal S(\alpha_t)\\M_t\in\mathcal S(\mu_t)}}\!\!\!\!\! F_{A_tB_t}^{M_t} \big) F_{C_fD_f}^{N_f} F_{AC_f}^{X_f} \big(\!\!\! \sum_{\substack{B_t\in\mathcal S(\beta_t)\\X_t\in\mathcal S(\xi_t)}}\!\!\!\!\! F_{B_tD}^{X_t} \big) \frac{a_\alpha a_{\beta_t}a_{\gamma_f}a_\delta}{a_\mu a_\nu}\\
&= \sum_{\substack{\alpha_t,\beta_t\\\gamma_f,\delta_f}}q^{\langle\xi_f,\xi_t\rangle-\langle\alpha,\delta\rangle} F_{\alpha_t\beta_t}^{\mu_t} F_{\gamma_f\delta_f}^{\nu_f} F_{\alpha\gamma_f}^{\xi_f} F_{\beta_t\delta}^{\xi_t} \frac{a_\alpha a_{\beta_t}a_{\gamma_f}a_\delta}{a_\mu a_\nu} \frac{n_{\mu_t}n_{\xi_t}}{n_\alpha n_{\beta_t} n_{\nu_t}}\\
&= \sum_{\alpha,\beta,\gamma,\delta}q^{\langle\xi_f,\xi_t\rangle-\langle\alpha,\delta\rangle} F_{\alpha\beta}^\mu F_{\gamma\delta}^\nu F_{\alpha\gamma}^{\xi_f} F_{\beta\delta}^{\xi_t} \frac{a_\alpha a_\beta a_\gamma a_\delta}{a_\mu a_\nu} \frac{n_\mu n_\xi}{n_\alpha n_\beta n_\nu}\\
&= F_{\mu\nu}^\xi n_\xi/n_\nu.
\end{align*}
Dually
$$\sum_{\substack{N\in\mathcal S(\nu)\\X\in\mathcal S(\xi)}}\!\!\!F_{MN}^X = F_{\mu\nu}^\xi n_\xi/n_\mu.$$
This proves the existence of the Hall polynomial $F_{\mu\nu}^\xi$. Note also that $n_\nu=n_\delta$ and $n_\gamma=1$, so that we can write $F_{\mu\nu}^\xi$ as in the statement of the proposition.

In the case that all indecomposable homogeneous regular modules are contained in $\mathcal F$, the proof goes through \textit{mutatis mutandis}.
\end{proof}

\begin{Lem}
It is enough to prove the existence of polynomials $F_{\mu\nu}^\xi$ when $\xi$ is indecomposable preprojective, indecomposable preinjective, non-homogeneous regular in a single tube, or homogeneous regular.
\end{Lem}

\begin{proof}
If $\xi$ can be expressed as $\xi_f\oplus\xi_t$ for some split torsion pair $(\mathcal T,\mathcal F)$ with both $\xi_f$ and $\xi_t$ non-zero and with all indecomposable homogeneous regular modules contained in either $\mathcal F$ or $\mathcal T$, then the previous proposition together with induction on dimension vector implies the result. In particular, we immediately reduce to the cases when $\xi$ is either preprojective, non-homogeneous regular in a single tube, homogeneous regular or preinjective.

If $\xi$ is decomposable preprojective, then we can find a section of the Auslander-Reiten quiver containing just one isotypic summand of $X\in\mathcal S(\xi)$, and with all other summands lying to the left. We can use this section to define a torsion pair, and so we can again apply our reduction.
Suppose that $\xi$ is isotypic preprojective, say $X=P^{r+1}\in\mathcal S(\xi)$ for some indecomposable preprojective $P$. Since $P$ is exceptional, we have (as in the proof of Theorem \ref{Ringel}) that
$$F_{MN}^{P^{r+1}} = \sum_{A,B,C,D}q^{r-\langle A,D\rangle}F_{AB}^MF_{CD}^NF_{AC}^PF_{BD}^{P^r}a_Aa_Ba_Ca_D/a_Ma_N.$$
We note that $N$, $C$ and $D$ are all preprojective, so discrete. Thus we can again use induction to deduce that
$$F_{MN}^{P^{r+1}} = \sum_{\alpha,\beta,\gamma,\delta}q^{r-\langle\alpha,\delta\rangle}F_{\alpha\beta}^\mu F_{\gamma\delta}^\nu F_{\alpha\gamma}^\pi F_{\beta\delta}^{\pi^r}\frac{a_\alpha a_\beta a_\gamma a_\delta}{a_\mu a_\nu}\frac{1}{n_\alpha n_\beta} \quad\textrm{for all }M\in\mathcal S(\mu).$$
Thus it is enough to assume that $\xi$ is indecomposable preprojective.
Similarly it is enough to consider indecomposable preinjective $\xi$.
\end{proof}
\begin{Lem}
It is enough to prove the existence of $F_{\mu\nu}^\xi$ when $\xi$ is regular and either homogeneous or else lying in a single non-homogeneous tube.
\end{Lem}
\begin{proof}
Suppose that $X\in\mathcal S(\xi)$ is indecomposable preprojective. Let $R=\mathrm{rad}(X)$ and $T=\mathrm{top}(X)$, so $R$ is preprojective and $T$ is semisimple, and both are discrete. If $R=0$, then $X$ is simple projective and the result is clear, so assume that both $R$ and $T$ are non-zero. Consider the left hand side of Green's Formula. We have that if $F_{TR}^E$ is non-zero, then either $E\cong X$ and $F_{TR}^X=1$, or else $E$ is decomposable and contains a preprojective summand (since it has the same defect as $X$). Thus
$$F_{MN}^X/a_X = \sum_EF_{MN}^EF_{TR}^E/a_E - \sum_{L\,\mathrm{decomp}} F_{MN}^LF_{TR}^L/a_L.$$
Conside the second sum. Since $L$ is decomposable and contains a preprojective summand, and since $R$ and $T$ are discrete, we can follow the proof of the previous lemma and apply induction to deduce the existence of a universal polynomial $F_{\tau\rho}^\lambda$ such that
$$F_{\tau\rho}^\lambda = F_{TR}^L \quad\textrm{for all }L\in\mathcal S(\lambda).$$
Similarly, since $N$ is preprojective, so discrete, there exists a universal polynomial $F_{\mu\nu}^\lambda$ such that
$$F_{\mu\nu}^\lambda n_\lambda/n_\mu = \sum_{L\in\mathcal S(\lambda)} F_{MN}^L \quad\textrm{for all }M\in\mathcal S(\mu).$$
Hence
$$\sum_{L\,\mathrm{decomp}}F_{MN}^LF_{TR}^L/a_L = \sum_{\lambda\neq\xi} F_{\mu\nu}^\lambda F_{\tau\rho}^\lambda n_\lambda/n_\mu a_\lambda.$$
Now consider the first sum. Using Green's Formula, we can rewrite it as
$$\sum_{A,B,C,D}q^{-\langle A,D\rangle}F_{AB}^MF_{CD}^NF_{AC}^TF_{BD}^R\frac{a_Aa_Ba_Ca_D}{a_Ma_Na_Ta_R}.$$
We note that $A$ and $C$ must both be semisimple and $D$ must be preprojective, so all three are discrete. Hence by induction we have universal polynomials such that
$$F_{BD}^R = F_{\beta\delta}^\rho/n_\beta, \quad F_{AC}^T=F_{\alpha\gamma}^\tau, \quad F_{CD}^N = F_{\gamma\delta}^\nu, \quad \sum_{B\in\mathcal S(\beta)}F_{AB}^M = F_{\alpha\beta}^\mu,$$
where the first holds for all $B\in\mathcal S(\beta)$ and the last holds for all $M\in\mathcal S(\mu)$.
Putting this together we obtain a universal polynomial such that $F_{MN}^X=F_{\mu\nu}^\xi$ for all $M\in\mathcal S(\mu)$.
If $\xi$ is indecomposable preinjective, then we can apply a similar argument using $\mathrm{soc}(X)$ and $X/\mathrm{soc}(X)$.
\end{proof}

\textbf{Remark.}
In fact, if $X$ lies in a single non-homogeneous tube, we may further assume that $X$ is indecomposable. For, $X$ is discrete and if $X=X_1\oplus X_2$ is decomposable, then every other extension $E$ of $X_2$ by $X_1$ lies in the same tube and satisfies $\lambda(E)<\lambda(X)$ (in the notation of Lemma \ref{cyclic-ext}). We can thus use induction on $\lambda(X)$ and Green's Formula as in the proof of Theorem \ref{cyclic}.

\section{Reductions using Associativity}\label{section:associativity}

In the previous section we reduced the problem of finding Hall polynomials to the special case of $F_{MN}^X$ for $X$ either homogeneous regular or non-homogeneous regular lying in a single tube.

We now wish to improve this to the case when $M$ is simple preinjective and $N$ is indecomposable preprojective. Unfortunately it seems difficult to do this by applying reflection functors, since we have often used induction on dimension vector. Instead we will use associativity together with Theorem \ref{Schofield}.

\begin{Lem}\label{indec-MN}
It is enough to consider $F_{\mu\nu}^\xi$ when $\nu$ is indecomposable preprojective, $\mu$ is indecomposable preinjective and $\xi$ is either homogeneous regular or else non-homogeneous regular in a single tube.
\end{Lem}
\begin{proof}
We showed in the last section that it is enough to consider the case when $\xi$ is homogeneous regular or non-homogeneous regular in a single tube.
We use induction on the dimension vector of $\xi$ and the defect of $\mu$. We note that if $\partial(\mu)=0$, then both $M$ and $N$ are regular. Thus either $M$, $N$ and $X$ are all homogeneous regular, so there exists a universal polynomial by Theorem \ref{gen_cyclic}, or else $M$, $N$ and $X$ all lie in a single non-homogeneous tube, so there exists a universal polynomial by Theorem \ref{cyclic}.
Consider $F_{MN}^X$ for $X\in\mathcal S(\xi)$ and assume that $\partial(\mu)>0$. Suppose further that we can find a split torsion pair such that $\mu_f$ and $\mu_t$ are both non-zero, and such that all homogeneous regulars are torsion-free. Note that
$$F_{MN}^X = \sum_L F_{M_fL}^XF_{M_tN}^L.$$
Since $\dimv L<\dimv X$ and $M_t$ is discrete, there exists a Hall polynomial $F_{\mu_t\nu}^\lambda$ such that
$$\sum_{N\in\mathcal S(\nu)}F_{M_tN}^L = F_{\mu_t\nu}^\lambda \quad\textrm{for all }L\in\mathcal S(\lambda).$$
Since $0\leq\partial(M_f)<\partial(M)$, there exists a Hall polynomial $F_{\mu_f\lambda}^\xi$ such that
$$\sum_{\substack{M_f\in\mathcal S(\mu_f)\\L\in\mathcal S(\lambda)}}F_{M_fL}^X = F_{\mu_f\lambda}^\xi \quad\textrm{for all }X\in\mathcal S(\xi).$$
Thus
$$F_{\mu\nu}^\xi := \sum_{\substack{M\in\mathcal S(\mu)\\N\in\mathcal S(\nu)}}F_{MN}^X = \sum_\lambda F_{\mu_f\lambda}^\xi F_{\mu_t\nu}^\lambda.$$
Similarly
$$\sum_{\substack{M\in\mathcal S(\mu)\\X\in\mathcal S(\xi)}}F_{MN}^X = F_{\mu\nu}^\xi n_\xi/n_\nu \quad\textrm{and}\quad
\sum_{\substack{N\in\mathcal S(\nu)\\X\in\mathcal S(\xi)}}F_{MN}^X = F_{\mu\nu}^\xi n_\xi/n_\mu.$$
This proves the existence of the Hall polynomial $F_{\mu\nu}^\xi$.

We have reduced to the case when $M$ is isotypic preinjective. Suppose that $M=I^{r+1}$ with $I$ indecomposable preinjective. Since $I$ is exceptional,
$$F_{II^r}^L = \delta_{LM}\frac{q^{r+1}-1}{q-1},$$
so associativity gives
$$F_{MN}^X = \frac{q-1}{q^{r+1}-1}\sum_L F_{IL}^XF_{I^rN}^L.$$
Again, $\dimv L<\dimv X$, so we have the Hall polynomial $F_{\iota^r\nu}^\lambda$. Also, $0\leq\partial(I)<\partial(M)$, so we have the Hall polynomial $F_{\iota\lambda}^\xi$, and hence the Hall polynomial $F_{\mu\nu}^\xi$.

The dual arguments clearly work for $N$, and since $0=\partial(X)=\partial(M)+\partial(N)$, we see that we can always reduce to the case when both $M$ is indecomposable preinjective and $N$ is indecomposable preprojective.
\end{proof}

The next lemma is a nice generalisation of what happens in the Kronecker case \cite{Szanto}.

\begin{Lem}\label{shift}
Let $X$ be either homogeneous regular or non-homogeneous regular lying in a single tube. If $M$ and $N$ are indecomposable with $\partial(M)=1$ and $\dimv M>\delta$, then
$$F_{MN}^X=F_{M'N'}^X,$$
where $M'$ and $N'$ are indecomposable and $\dimv M-\dimv M'=\delta=\dimv N'-\dimv N$.
\end{Lem}

\begin{proof}
By assumption, there exists an embedding of the module category of the Kronecker quiver whose image contains $M$. (For, take an indecomposable preinjective module $I$ such that $\dimv I<\delta$ and $\dimv M-\dimv I\in\mathbb Z\delta$. Then there exists an orthogonal exceptional pair $(I,P)$ with $[I,P]^1=2$.) In this way, we see that there exists a short exact sequence of the form
$$0\to R\to M\to M'\to 0$$
for any indecomposable regular $R$ such that $\dimv R=\delta$ and $[M,T]^1\neq0$, where $T$ is the regular top of $R$; equivalently $[S,M]\neq0$, where $S$ is the regular socle of $R$. This occurs for precisely one such indecomposable in each non-homogeneous tube. It is easily seen from Riedtmann's Formula that
$$F_{M'R}^M=1, \quad F_{M'R}^{R\oplus M'}=q \quad\textrm{and}\quad F_{RM'}^{R\oplus M'}=1.$$
Dually, we have a short exact sequence
$$0\to N\to N'\to R'\to 0$$
for any indecomposable regular $R'$ such that $[S',N]^1\neq0$, where $S'$ is the regular socle of $R'$.

In particular we can take $R=R'$ for any homogeneous regular indecomposable $R$ of dimension vector $\delta$. It may be, however, that there are no such homogeneous regular modules of dimension vector $\delta$. For example, take $Q$ a quiver of type $\widetilde{\mathbb D}_4$ and $k$ a field with two elements.

Suppose, therefore, that there are no homogeneous regular modules of dimension vector $\delta$. By assumption we have a short exact sequence
$$0\to N\to X\to M\to 0$$
with $X$ either homogeneous regular or non-homogeneous regular lying in a single tube. Take a regular simple $S$ in a non-homogeneous tube different from that containing $X$, and apply $\Hom(S,-)$. Since there are no homomorphisms or extensions between regular modules in distinct tubes, we see that $[S,M]=[S,N]^1$. Take $S$ such that $[S,M]=1$ and set $R$ to be the indecomposable regular module of dimension vector $\delta$ and with regular socle $S$.

In all cases, we have found a module $R$ for which there exist exact sequences
$$0\to R\to M\to M'\to 0 \quad\textrm{and}\quad 0\to N\to N'\to R\to 0.$$
Associativity of Hall numbers now gives
$$F_{MN}^X+qF_{R\oplus M'N}^X = \sum_LF_{M'R}^LF_{LN}^X = \sum_LF_{M'L}^XF_{RN}^L = F_{M'N'}^X+qF_{M'N\oplus R}^X.$$
Similarly
$$F_{R\oplus M'N}^X = \sum_LF_{RM'}^LF_{LN}^X = \sum_LF_{RL}^XF_{M'N}^L$$
and
$$F_{M'N\oplus R}^X = \sum_LF_{M'L}^XF_{NR}^L = \sum_LF_{M'N}^LF_{LR}^X.$$
Thus
$$F_{MN}^X-F_{M'N'}^X = q\sum_LF_{M'N}^L\big(F_{LR}^X-F_{RL}^X\big).$$
If $R$ is homogeneous regular, then it is well-known (using the natural duality on representations) that the Hall numbers are `symmetric'; that is, $F_{LR}^X=F_{RL}^X$. On the other hand, if $R$ is non-homogeneous and lying in a different tube from $X$, then $F_{LR}^X=F_{RL}^X=0$. In all cases we get that
$$F_{MN}^X=F_{M'N'}^X$$
as required.
\end{proof}

\begin{Lem}
It is enough to prove the existence of Hall polynomials in the case when $\mu$ is simple preinjective, $\nu$ is indecomposable preprojective and $\xi$ is regular and either homogeneous or contained in a single non-homogeneous tube.
\end{Lem}

\begin{proof}
By Lemma \ref{shift}, we may assume that $M$ is indecomposable preinjective and either $\partial(M)\geq2$ or else $\partial(M)=1$ and $\dimv M<\delta$. If $M$ is not simple, then Theorem \ref{Schofield} yields an orthogonal exceptional pair $(A,B)$ and a short exact sequence of the form
$$0\to B^b\to M\to A^a\to 0.$$
If $d=\dim\Ext^1(A,B)=2$, then by Lemma \ref{Lem_Schofield} and using that $M$ is indecomposable preinjective, we have $\dimv M=r\delta+\dimv A>\delta$ where $r=b=a-1$, and $\partial(M)=\partial(A)=1$, a contradiction. Thus $d=1$ and, since there are only three indecomposable modules in $\mathcal F(A,B)$, $a=b=1$. It is easily seen that $F_{AB}^L=1$ for both $L=M$ and $L=A\oplus B$. Hence
$$F_{MN}^X+F_{A\oplus BN}^X = \sum_LF_{AB}^LF_{LN}^X = \sum_LF_{AL}^XF_{BN}^L.$$
Note also that
$$F_{A\oplus BN}^X=\sum_LF_{BA}^LF_{LN}^X=\sum_LF_{BL}^XF_{AN}^L.$$
Putting this together we obtain
$$F_{MN}^X = \sum_LF_{AL}^XF_{BN}^L - \sum_LF_{BL}^XF_{AN}^L.$$
Now, since in both sums $\dimv L<\dimv X$, we have Hall polynomials $F_{\beta\nu}^\lambda$ and $F_{\alpha\nu}^\lambda$. Since $M\succ A$, we can use induction on the partial order for indecomposable preinjective modules to deduce the existence of the Hall polynomial $F_{\alpha\lambda}^\xi$. Finally, there are three possibilities for $B$: if $B$ is preprojective, then $F_{BL}^X=0$; if $B$ is (non-homogeneous) regular, then there is a Hall polynomial for $F_{BL}^X$ since all three modules are regular; if $B$ is preinjective, then $\partial(M)>\partial(A),\partial(B)>0$ and we can use induction on the defect to obtain a Hall polynomial $F_{\beta\lambda}^\xi$.

We have shown that, for all $X\in S(\xi)$,
$$F_{MN}^X=\sum_\lambda F_{\alpha\lambda}^\xi F_{\beta\nu}^\lambda - \sum_\lambda F_{\beta\lambda}^\xi F_{\alpha\nu}^\lambda,$$
and hence the Hall polynomial $F_{\mu\nu}^\xi$ exists.
\end{proof}

We now show that Hall polynomials exist when $M$ is simple preinjective, $N$ is indecomposable preprojective, and $X$ is regular and either homogeneous or contained in a single non-homogeneous tube, thus completing the proof of the Main Theorem.

Let $d=[N,M]$. Then there is an epimorphism $N\twoheadrightarrow M^d$, unique up to an automorphism of $M^d$. Let $P$ be the kernel of this map. We thus have a short exact sequence
$$0\to P\to N\to M^d\to 0$$
and $F_{M^dP}^N=1$. Applying $\Hom(-,M)$ yields $[P,M]=0$, whereas $\Hom(-,N)$ yields $[P,N]^1=0$. Finally, using $\Hom(P,-)$, we obtain that $[P,P]^1=0$ and hence that $P$ is rigid. Consider now the push-out diagram
$$\begin{CD}
&& P @= P\\
&& @VVV @VVV\\
0 @>>> N @>>> X @>>> M @>>> 0\\
&& @VVV @VVV @|\\
0 @>>> M^d @>>> M^{d+1} @>>> M @>>> 0
\end{CD}$$
where the bottom row is split since $M$ is exceptional. We note that $[X,M]=d+1$ and $F_{M^{d+1}P}^X=1$.

By associativity,
$$F_{MN}^X + \sum_{L\,\mathrm{decomp}}F_{ML}^XF_{M^dP}^L = \sum_L F_{ML}^XF_{M^dP}^L = \sum_L F_{MM^d}^LF_{LP}^X = \frac{q^{d+1}-1}{q-1},$$
so that
$$F_{MN}^X = \frac{q^{d+1}-1}{q-1} - \sum_{L\,\mathrm{decomp}}F_{ML}^XF_{M^dP}^L.$$

Since $\dimv L<\dimv X$ and $M$ and $P$ are discrete, there exists a Hall polynomial $F_{\mu^d\pi}^\lambda$ such that
$$F_{M^dP}^L = F_{\mu^d\pi}^\lambda \quad\textrm{for all }L\in\mathcal S(\lambda).$$
Since $L$ is decomposable and has a preprojective summand ($\partial(L)=\partial(N)<0$), Lemma \ref{indec-MN} gives us the Hall polynomial $F_{\mu\lambda}^\xi$. Thus there exists a Hall polynomial $F_{\mu\nu}^\xi$ as required.

This completes the proof of the Main Theorem.

\section{Wild quivers}

We finish by offering a definition of Hall polynomials for wild quivers.

Let $Q$ be an arbitrary quiver. We first need a combinatorial partition of the set of isomorphism classes of $kQ$-modules, by which we mean some combinatorial `properties' of modules such that the sets $\mathcal S(\alpha,k)$ yield a partition of the set of isomorphism classes of modules, where $\mathcal S(\alpha,k)$ consists of those $kQ$-modules having property $\alpha$. Moreover, we need universal polynomials $a_\alpha$ and $n_\alpha$ for each $\alpha$ such that, for any finite field $k$ with $|k|=q$,
$$a_\alpha(q)=|\Aut(A)| \quad\textrm{for each }A\in\mathcal S(\alpha,k) \quad\textrm{and}\quad n_\alpha(q)=|\mathcal S(\alpha,k)|.$$
Finally, we would like that there are only finitely many such classes of a given dimension vector.

Some further properties which would clearly be of interest are that, over an algebraically closed field $k$, each $\mathcal S(\alpha,k)$ is a locally closed, smooth, irreducible, $\mathrm{GL}$-invariant subvariety of the representation variety, that each $\mathcal S(\alpha,k)$ admits a smooth, rational, geometric quotient (in which case $n_\alpha$ counts the number of rational points of this quotient) and that the $\mathcal S(\alpha,k)$ yield a stratification of the representation varieties.

We say that Hall polynomials exist with respect to such a family if the following conditions are satisfied. Given $\alpha$, $\beta$ and $\gamma$, there exists a universal polynomial $F_{\alpha\beta}^\gamma$ such that
$$F_{\alpha\beta}^\gamma(q) = \sum_{\substack{A\in\mathcal S(\alpha,k)\\B\in\mathcal S(\beta,k)}}F_{AB}^C \quad\textrm{for all }C\in\mathcal S(\gamma,k)$$
and further that
\begin{align*}
n_\gamma(q) F_{\alpha\beta}^\gamma(q) &= n_\alpha(q) \sum_{\substack{B\in\mathcal S(\beta,k)\\C\in\mathcal S(\gamma,k)}} F_{AB}^C \quad\textrm{for all }A\in\mathcal S(\alpha,k)\\
&= n_\beta(q) \sum_{\substack{A\in\mathcal S(\alpha,k)\\C\in\mathcal S(\gamma,k)}} F_{AB}^C \quad\textrm{for all }B\in\mathcal S(\beta,k).
\end{align*}

We note that the results of Section \ref{section:Green} can be partially extended to such a situation, whereas the results of Section \ref{section:associativity} are particular to affine quivers since they require the notion of defect. This allows one to simplify the situation somewhat, but we are still left with considering Hall numbers involving at least two regular modules.

\end{document}